\documentclass[11pt,reqno]{amsart} 

\usepackage{amsmath}
\usepackage{amsthm}
\usepackage{amssymb}
\usepackage{amsfonts}
\usepackage{latexsym}
\usepackage{mathrsfs}
\usepackage{multicol}
\usepackage{fancyvrb}
\usepackage{lipsum}
\usepackage{latexsym,amsbsy,amsxtra,amscd,verbatim}

\usepackage{lmodern}
\usepackage[T1]{fontenc}

\usepackage[left=1cm, right=1cm]{geometry}
\usepackage{paralist}

\usepackage{color}
\usepackage[colorlinks,pagebackref=true]{hyperref}

\makeatletter

\def\@citecolor{blue}
\def\@linkcolor{blue}
\def\@urlcolor{blue}

\makeatother

\usepackage{multicol}
\usepackage{fancyvrb}
\usepackage[colorlinks]{hyperref}

\def\adj{\operatorname{adj}}

\newcommand{\core}{\operatorname{core}}

\def\C{\mathbb C}

\def\N{\mathbb N}

\def\m{\mathfrak m}

\def\@citecolor{blue}
\def\@linkcolor{blue}
\def\@urlcolor{blue}
\newcommand{\ov}{\overline}

\theoremstyle{plain}
\newtheorem{theorem}{Theorem}[section]
\newtheorem{proposition}[theorem]{Proposition}
\newtheorem{lemma}[theorem]{Lemma}
\newtheorem{corollary}[theorem]{Corollary}
\newtheorem{conjecture}[theorem]{Conjecture}
\newtheorem{principle}[theorem]{Principle}

\theoremstyle{definition}

\newtheorem{remark}[theorem]{Remark}
\newtheorem{example}[theorem]{Example}

\newcommand{\ncom}{\newcommand}
\ncom{\bib}{\bibitem}

\ncom{\limns}{\underset{ s \longrightarrow \infty }{\lim}}
\ncom{\limnr}{\underset{ r \longrightarrow \infty }{\lim}}
\ncom{\beqn}{\begin{eqnarray*}}
\ncom{\eeqn}{\end{eqnarray*}}
\ncom{\beq}{\begin{eqnarray}}
\ncom{\eeq}{\end{eqnarray}}
\ncom{\f}{\frac}
\ncom{\been}{\begin{enumerate}}
\ncom{\eeen}{\end{enumerate}}
\ncom{\olin}{\overline}
\ncom{\ulin}{\underline}

\title[]{
On the core and adjoint of the product of complete ideals in two-dimensional regular local rings}
\author[Clare D'Cruz]{Clare D'Cruz$^*$}
\address{Chennai Mathematical Institute, Plot H1 SIPCOT IT Park, Siruseri, 
Kelambakkam 603103, Tamil Nadu, 
India
} 
\email{clare@cmi.ac.in} 
\thanks{$*$ partially supported by a grant from Infosys Foundation and SERB  MATRICS grant  MTR/2023/000661.
}
\author[Saipriya Dubey] {Saipriya Dubey$^\dag$}
\address{Department of Mathematics, Indian Institute of Technology Dharwad, Dharwad, 580001, India}
\email{saipriyadubey@iitdh.ac.in}
\curraddr{Chennai Mathematical Institute, Plot H1 SIPCOT IT Park, Siruseri, 
Kelambakkam 603103, Tamil Nadu, 
India
}
\email{saipriyad@cmi.ac.in}

\author{Jugal K. Verma$^\dag$}
\address{Department of Mathematics, Indian Institute of Technology Bombay, Mumbai, 400076, India}
\email{jkv@iitb.ac.in}
\curraddr{Mathematics Department, 
IIT Gandhinagar, Palaj, Gujarat 382055, India}
\email{jugal.verma@iitgn.ac.in}
\thanks{$\dag$ partially supported by SPARC grant, Ministry of Human Resource and Development, Govt. of India.
}
\dedicatory{Dedicated to  Sudhir Ghorpade on the occasion of his $60$-th birthday}

\begin{document}

\begin{abstract}
Using joint reductions of complete ideals, we find  expressions for the core and adjoints of the product of complete ideals in a two-dimensional regular local ring. We also compute  their colengths. Our results strengthen a generalization of the Brian\c{c}on-Skoda theorem due to D.~Rees and J.~D.~Sally.
\end{abstract}

\maketitle
\section{Introduction}
Throughout this paper,  we will assume that $(R, \m)$ is a  Noetherian local ring of dimension $d$ with infinite residue field $R/ \m$. 
 Let  $I$ be an ideal minimally generated by $\mu$  elements and  let $\ov{I}$ denote the integral closure of  $I.$ 
One of the weakest forms of the  Brian\c{c}on-Skoda theorem states that  if $(R,\m)$ is a regular local ring, then  $\ov{I^{n+\mu-1}}\subseteq I^n$ for all $n\geq 1$ (\cite{BS1974}, \cite{LS1981}, \cite{LT1981}). 
If   $I$ is an  $\m$-primary ideal,  then any  minimal reduction $J$ of $I$ is generated by  $d$ elements and hence  $\ov{I^d}=\ov{J^d} \subseteq J$. It follows that the intersection of all minimal reductions of $I$, called the core of $I$, denoted by $\core (I)$,  contains $\ov{I^{d}}$. This implies that  for an $\m$-primary ideal $I$,   $\core (I)$   is also $\m$-primary.  The core of an ideal  was introduced by J.~D.Sally and  was defined in a paper of D.~Rees and J.~D.~Sally in \cite{rees-sally}. They showed that if $\{x_1, x_2, \ldots, x_d\}$ is a joint reduction of a  set of  $\m$-primary ideals $\{I_1, I_2, \ldots, I_d\}$,  
 then 
$\ov{I_1 I_2 \cdots I_d}\subseteq (x_1, x_2, \ldots, x_d)$ \cite{rees-sally}. Their result is an  improvement of  the Brian\c{c}on-Skoda theorem for a family of ideals. 
An interesting consequence of their result  is that  the intersection of all joint reductions of a set of  $d$ $\m$-primary ideals  is also an $\m$-primary ideal.

In    \cite{lipman-adjoint}, J.~Lipman used  adjoints of ideals   to give  a generalization of the Brian\c{c}on-Skoda theorem. 
Let $(R, \m)$ be a regular local ring and let  $K$ be the fraction field of $R$.
The adjoint of an ideal $I$ is the ideal
$$
\adj(I)=\bigcap_{v}\{r\in K\mid v(r)\geq v(I)-v(J_{R_v/R})\},
$$
 where the intersection is taken over all prime divisors   $v$  of $R$,  and  the Jacobian ideal $J_{R_v/R}$ is the  $0$-th Fitting ideal of the $R_v$-module of K\"{a}hler diﬀerentials $\Omega^1_{R_v/R}$. 
 The adjoint of $I$ is  a complete ideal  and  $I \subseteq  \ov{I}\subseteq \adj(I)$ \cite[Remark~1.2(b)]{lipman-adjoint}. 
J.~Lipman    proposed the following conjecture. 
\begin{conjecture}\label{adjointconjecture}
\cite[Conjecture~1.6]{lipman-adjoint}
Let $I$ be an ideal in a regular local ring $R$ of analytic spread  $\ell(I)$.   Then  for all $ n \geq \ell(I)-1$, 
$$
\adj(I^{n+1})=I \cdot \adj(I^{n}).
$$
\end{conjecture}
A positive answer to this conjecture would give  quick proofs for  the  Brian\c{c}on-Skoda  type theorems. J.~Lipman  proved that Conjecture~\ref{adjointconjecture} holds true  for ideals in a   two-dimensional regular local ring (see \cite[(2.3)]{lipman-adjoint}).  
In  \cite{hs06},  C.~Huneke and I.~Swanson proved an interesting result which relates the core of an ideal  and adjoint of an ideal. 
They showed that if $I$ is an $\m$-primary  complete ideal   in a regular local ring of dimension two, then  $\core(I) = I \cdot \adj(I)$ 
\cite[Theorem~3.14]{hs06}. 
It follows that   if $J$ is a minimal reduction of $I$, then 
 $$
 \ov{I^2}=\ov{J^2}\subseteq \adj(J^2)=\adj(I^2)=J \cdot \adj(J) \subseteq \core(I).
 $$
This gives a stronger version of Brian\c{c}on-Skoda theorem, namely $\adj(I^2) \subseteq \core(I).$
  
 To prove Conjecture~\ref{adjointconjecture},  J.~Lipman formulated a conjecture  which he called  the `Vanishing Conjecture' \cite[Conjecture~2.2]{lipman-adjoint}. He also remarked  that the `Vanishing  Conjecture'  implies Conjecture~\ref{adjointconjecture}. 
 The `Vanishing Conjecture' was solved by S.~D.~Cutkosky \cite[Theorem A3]{lipman-adjoint} for  local domains that are of essentially finite type over a field of characteristic zero and hence in this setting,   
  Conjecture~\ref{adjointconjecture} also holds true.  In \cite{hs06}, several interesting properties of adjoint and core were proved for $\m$-primary ideals in a two-dimensional regular local ring. Their proofs are simple and give an explicit  description of the adjoint and core. 
  In  \cite{hubl-swanson}, R. H\"{u}bl and I. Swanson proved that the adjoint of a  generalized monomial ideal is a generalized  monomial ideal. 
      In \cite{KM},   M.~Kummini and S.~Masuti proved that if $R$ is a  regular local ring of dimension three and  if  $X=\text{Proj} (\oplus_{n=0}^\infty \ov{I^n})$ is pseudo-rational,  then $\adj{(I^n)}=I \cdot\adj{(I^{n-1})}$ for all $n\geq 3.$ 

Motivated by earlier results on the adjoint,  one can ask  the following question: Is it possible to express  $\adj(IJ)$ in terms of $\adj(I)$ and $\adj(J)$? 
In Section~3 we address this question. In Section~4 we express  $\core(IJ)$ in terms of $\core(I)$ and $\core(J)$. 
The main results in this paper are:
\begin{theorem}
Let  $I$ and $J$ be $\m$-primary  ideals in a  two-dimensional regular local ring $(R,\m)$.  Let $\{a,b\}$ be  a joint reduction of $\{I, J\}$. Then for all $r,s \geq 1,$
 \begin{align}
 \label{section one main adj}
 \adj(I^rJ^s)
 =a^r \cdot \adj(J^s) + b^s \cdot \adj(I^r) = I^r \cdot \adj (J^s) + \adj (I^r)  \cdot J^s. 
 \end{align}
 \end{theorem}

\begin{theorem}
\label{joint-reduction-zero-core-intro}
Let $I$ and $J$ be $\m$-primary complete ideals in a two-dimensional  regular local ring $(R, \m)$.  
Let $\{a,b\}$ be a joint reduction of $\{I,J \}$. Then  for all $r,s \geq 1$, 
 \begin{align}
  \label{section one main core-intro}
 \core(I^rJ^s) = a^{2r} \core(J) + b^{2s} \core(I) = I^{2r} \core(J^s) + \core(I^r) J^{2s}.
 \end{align}
 \end{theorem}
By Lemma~\ref{lemma-containment} (resp. Lemma~\ref{lemma-containment-core}) one can show that  $I^r \cdot \adj (J^s) + \adj (I^r)  \cdot J^s \subseteq \adj(I^r J^s)$ 
(resp.  $I^{2r} \core(J^s) + \core(I^r) J^{2s} \subseteq  \core(I^rJ^s)$). However, equality  is not obvious and we use 
powerful result by Hoskin and Deligne  which computes the colength of an $\m$-primary  complete ideal in a two-dimensional regular local  ring (see \cite[Theorem 2.13]{Del73}, \cite[Theorem 5.1]{Hos56}). Several researchers have reproved the result of Hoskin and Deligne for colength of complete ideals (\cite[14.5.4]{huneke-swanson},   \cite[Theorem 3.10]{john-verma},  \cite{kodiyalam}  and \cite[Theorem 3.1]{lipman-fs}). 
 
 For any two ideals $I,J$ in a regular local ring $R$,  we say that the   subadditivity property holds for adjoints of ideals if $\adj(IJ ) \subseteq  \adj(I) \cdot \adj(J )$. 
In  \cite{hubl-swanson}   R. H\"{u}bl and I. Swanson proved  the subadditivity property of adjoints of ideals  in some special cases.
In \cite{tak-wat},  S~Takagi and K.~Watanabe proved the subadditivity property for multiplier ideals and thus deduced the subadditivity property for adjoints  of  ideals in a  two-dimensional regular local ring.
We would like to remark that from (\ref{section one main adj}) and   Lemma~\ref{lemma-containment} it follows that if $I$ and $J$ are $\m$-primary ideals in a two-dimensional regular local ring, then $\adj(IJ) \subseteq \adj(I) \cdot \adj(J)$, thus giving an alternate proof 
 for the subadditivity property of adjoints of $\m$-primary ideals in a two-dimensional regular local ring. 
However,  the subadditivity property for adjoints of  ideals  is not known in general. 
From (\ref{section one main core-intro}) and Lemma~\ref{lemma-containment-core}  it follows that if $I$ and $J$ are complete ideals in a two-dimensional regular local ring, then  $\core(IJ) \subseteq \core(I) \cdot \core(J)$. 
This proves the subadditivity of  property of the core of $\m$-primary complete ideals in a two-dimensional regular local ring.  This is a new result which was not known earlier.

In Section 5, we illustrate our results with examples. 

We thank the referee for the careful reading of our manuscript and the valuable suggestions. 

\section{Preliminaries}
In this section, we state several basic definitions and results  that will be useful in the later sections. 
For all undefined terms,  we request  the reader to refer  to   W.~Bruns and J.~Herzog \cite{bruns-herzog} and H. Matsumura \cite{matsumura} . 

\subsection{Reductions and joint reductions}
 Let  $I$ be an ideal in a local ring $(R,\m)$ of dimension $d$. An ideal $J$ is called a reduction of  $I$  if $J\subseteq I$ and  there exists $n \geq 0$ such that $J I^n = I^{n+1}$. 
Any minimal reduction of $I$ is generated by $\ell(I)$ elements. Recall that  $\ell(I)$, the analytic spread of $I$,   is  the Krull dimension  of the fiber ring 
 $F(I)=\bigoplus_{n=0}^\infty I^n/\m I^n.$ 
 The analytic spread of an ideal $I$  is bounded below by the height of $I$ and  above by $d$. 
   Hence, if $I$ is $\m$-primary, then $\ell(I)=d.$
   
 An element $r \in R$ is said to be  integral over $I$ if there exists an integer $n \in \N$ and elements $a_i \in  I^i$,    $i = 1,2,\ldots, n$ such that 
$$
r^n + a_1 r^{n-1}  + a_2 r^{n-2} +\cdots  + a_{n-1} r+ a_n =0.
$$ 
 The  set of all elements that are integral over $I$ is called the integral closure of $I$ and is denoted by $\olin{I}$.  The ideal $I$ is called a complete ideal if 
 $I=\olin{I}$.  An  ideal  $J \subseteq I$ is a reduction of $I$ if and only if $I\subseteq \olin{J}.$ 

Let $I$ be an $\m$-primary ideal of  $(R,\m).$ The Hilbert function of $I$ is defined as $H_I(n):=\lambda(R/I^n).$ For all large $n$,  this function  is given by  a polynomial  $P_I(n) $,  called the Hilbert polynomial of $I$ and    can written as
\[
   P_I(n)
=e(I)\binom{n+d-1}{d}-e_1(I)\binom{n+d-2}{d-1}+\cdots+(-1)^d e_d(I).
\]
The integers  $e(I), e_1(I), \ldots, e_d(I),$ are called the Hilbert coefficients of $I.$  If $J\subseteq I$ is a reduction of $I$, then $e(I)=e(J). $ D.~Rees proved that if $R$ is formally equidimensional  and  
$J \subseteq I$ are $\m$-primary ideals then  $e(I)=e(J)$ if and only if  $\ov{I} =  \ov{J}$  \cite[Theorem~3.2]{Rees1961}.
In  \cite{Ree84},  D.~Rees introduced joint reductions of ideals to study  mixed multiplicities of ideals. 

Let  $I_1, I_2, \ldots, I_g$ be $\mathfrak{m}$-primary ideals in $(R, \m)$. The Hilbert function of $\underline{I}=I_1, I_2, \ldots, I_g$ is the function 
$$
    H_{\underline{I}}(n_1, n_2, \ldots, n_g)
= \lambda \left( \frac{R}{I_1^{n_1}I_2^{n_2}\cdots I_g^{n_g}} \right).
$$
For all large $n_1, n_2, \ldots, n_g,$ $H_{\underline{I}}(n_1, n_2, \ldots, n_g)$ is given by a polynomial  $P_{\underline{I}}(n_1, n_2, \ldots, n_g)$ which can be  written as
$$
   P_{\underline{I}}(n_1, n_2, \ldots, n_g)
=\sum_{j_1+j_2+\dots+j_g\leq d} e(j_1,j_2, \ldots,  j_g)
   \binom{n_1+j_1-1}{j_1}\binom{n_2+j_2-1}{j_2}\cdots \binom{n_g+j_g-1}{j_g}
   $$
where  $e(j_1, j_2, \ldots , j_g)$ are integers.  The integers $e(j_1, j_2, \ldots, j_g)$ for which $j_1+j_2+\cdots+j_g=d$ are called the mixed multiplicities of the set of ideals  $\{I_1, I_2, \ldots, I_g\}$.  In \cite{Ree84}, D.~Rees introduced joint reductions  and expressed mixed multiplicities in terms of joint reductions. 
A set of elements $\{x_1, x_2, \ldots, x_d\}$ is called a joint reduction of the set of ideals $\{I_1, I_2, \ldots, I_d\}$ if $x_j\in I_j$ for $j=1, 2, \ldots, d$ and the ideal 
$
{\displaystyle \sum_{j=1}^d I_1 I_2\cdots I_{j-1}x_jI_{j+1}\cdots I_d}
$
is a reduction of $I_1I_2\cdots I_d.$ 
D.~Rees proved that if $\{x_1, x_2, \ldots, x_d\}$ is a joint reduction of  $\{I_1, I_2, \ldots, I_d\}$, then $e(1,1,\ldots, 1)=e(x_1, x_2, \ldots, x_d)$ \cite[Theorem~2.4]{Ree84}.  The converse of Rees' result was proved by I.~Swanson  for formally equidimensional  local rings  (\cite[Theorem~3.7]{swanson},  \cite[Theorem 17.6.1]{huneke-swanson}). 

We state a well known result of  D.~Rees which will be  frequently used in the paper. 
\begin{lemma}
\label{rees-1961}
\cite[Lemma~2.4]{Rees1961}, \cite[page 398]{Ree84}
Let $(R, \m)$ be a Cohen-Macaulay local ring of dimension two and let $I,J$ be $\m$-primary ideals in $R$. Then 
$e(2,0)= e(I)$, $e(0,2) = e(J)$ and $e_1(I|J):= e(1,1) = (e(IJ) - e(I) - e(J) )/2 $.
\end{lemma}

\subsection{Complete ideals in  two-dimensional regular local rings} 
For the rest of this section we will assume that  $(R, \m)$ is a  two-dimensional regular local ring  and  $\m = (x, y)$. Any first local quadratic transform of $R$ is a localization of either $R[\m/x]$ or $R[\m/y]$ at a maximal ideal. 
 We say that a regular local ring $ (T, \m_{_{T}})$ birationally dominates $(R, \m)$ if  $R \subset T$, they have the same  fraction field  
and $\m_{_T} \cap  R = \m $. If $T$ birationally dominates $R$, we write it as $R \preceq T$. 
Any first local quadratic transform $(S, \m_{_S})$ of $(R, \m)$ is a  regular local ring  of dimension two that birationally dominates $R$. Any $n$-th local quadratic transform of $R$ is a first local quadratic transform of an $(n-1)$-st local quadratic transform of $R$.
 
The following remarkable result of S. Abhyankar \cite{abhyankar} plays an important role in the theory of complete ideals in two-dimensional regular local rings.
\begin{theorem}[Abhyankar, 1956]
Let $(R,\m_R) \prec (T,\m_{_{T}})$ be two-dimensional regular local rings  with  the same field of fractions $K$. Then there exists a unique sequence of two-dimensional regular local rings 
$$R=R_0  \prec R_1  \prec R_2 \prec \cdots \prec R_n=T\subset K,$$  where each 
 $R_i$  is a first local quadratic transform of $R_{i-1}$.
 \end{theorem}
Let $a$ be a nonzero element in  $R$. The $\m$-adic order of $a$, denoted by $o(a),$ is the largest power  $r$ such that $a \in \m^r$. 
The $\m$-adic order of an ideal $I$  denoted by $o(I)$,  is the largest power $r$ so that $I \subseteq \m^r$. 
  Let $x$ be a minimal generator of $\m = (x,y)$ and $T =R[\m/x]_{(x, y/x)}$. 
 For any ideal  $I$  in $R$,   the transform of $I$ in $T$ is defined as  $I^T:=x^{-r}(IT),$ where $r=o(I)$.  
We  put  $r_{_{T}}(I) := o(I^{^{T}}).$ 
  The point basis  of $I,$ denoted by $B(I)$,  is the set $\{r_{_T}(I) \mid R\preceq T\}$  \cite{lipman-fs}. 
 The residue field of $T$ is a finite algebraic extension of the residue field of $R$, and we denote it  by $[T:R]$. 

 We state a formula due to Hoskin and Deligne for the colength of an $\m$-primary complete ideal $I$ in  $R$ (\cite[Theorem 2.13]{Del73}, \cite[Theorem 5.1]{Hos56},   \cite[Theorem 3.10]{john-verma},  and \cite[Theorem 3.3]{lipman-vanishing}).
 We  also recall a   few well known results  which will be  used in later sections.
\begin{theorem}
\label{hoskin-deligne}
 Let $I,J$ be a complete $\m$-primary ideals of $R$. Then 
\been[\em (a)]
\item 
\label{hoskin-deligne-0}
\cite[Appendix~5, Theorem~$2^{\prime}$]{zar-sam} For all $r,s \geq 1$,  $I^rJ^s$ is a complete ideal. 

\item 
\label{hoskin-deligne-a}
{\em [Hoskin, Deligne]} For all  $n \geq 1$,  ${\displaystyle
\lambda \left( \frac{R}{I^n} \right) 
= \sum_{R\preceq T}\binom{r_{_T}(I^n)+1}{2} [T:R].
}$

\item 
\label{hoskin-deligne-b}
\cite[Theorem~3.7, Lemma~3.9]{john-verma}
${\displaystyle e(I) = \sum_{R\preceq T}(r_{_T}(I))^2 [T:R]}$ 
and ${\displaystyle e_1(I) = \sum_{R\preceq T} \binom{r_{_T}(I)}{2} [T:R]} = e(I)- \lambda (R/I)$.

\item
\label{hoskin-deligne-c}
\cite[Remark~2.6]{lipman-fs} 
${\displaystyle e_1(I|J)=\sum_{R\preceq T} r_{_T}(I)  r_{_T}(J)  }[T:R]$.
\eeen
\end{theorem}

\begin{theorem}
\label{lipman-length-twoideals}
Let $I$ and $J$ be $\m$-primary complete ideals in  $R$. Then 
\been[\em (a)]
\item 
\label{lipman-length-twoideals-1}
{\em (}\cite[Corollary~3.7]{lipman-fs}, \cite[Theorem~2.1, Theorem~3.2]{verma-nagoya}{\em)}
 For all $r,s \geq 0$
\beqn
\lambda \left(  \frac{R}{I^r J^s}\right)
= \lambda \left(  \frac{R}{I^r}\right) + \lambda \left(  \frac{R}{J^s}\right) + rse_1(I|J)
\eeqn

\item
\label{lipman-length-twoideals-3}
If $\{a,b\}$ is a joint reduction of $\{I,J\}$, then $e_1(I^r| J^s) = e(a^r, b^s) = rse(a, b)$. 

\item
\label{lipman-length-twoideals-4}
\cite[Lemma~2.5]{Ree84}  Let $I$, $J$ and $K$ be $\m$-primary ideals. Then $e_1(IJ|K) = e_1(I|K) + e_1(J|K)$.

\item
\label{lipman-length-twoideals-2}
For all $r,s \geq 1$, $e_1 (I^r|J^s) = rs e_1(I|J)$.

\item
\label{lipman-length-twoideals-5}
$e_1(\ov{I} | \ov{J}) = e_1(I |J)$.
\eeen
\end{theorem}
\begin{proof} 

(\ref{lipman-length-twoideals-3})  Since  $\{a,b\}$ is a joint reduction of $\{I, J\}$,  $\{a^r, b^s\}$ is a joint reduction of $\{I^r, J^s\}$. Hence applying \cite[Theorem~2.4(ii)]{Ree84} we get
$e_1(I^r | J^s) = e(a^r, b^s)$. From \cite[page 314]{lech} we get $e(a^r, b^s)= rs e(a,b)$.

(\ref{lipman-length-twoideals-2})  follows from (\ref{lipman-length-twoideals-4}).

(\ref{lipman-length-twoideals-5}) By Lemma~\ref{rees-1961},
\beqn
e_1( \ov{I} | \ov{J})
= \frac{1}{2} \left[ e( \ov{I} \ov{J})  -e(\ov{I}) -e(\ov{J})\right]
= \frac{1}{2} \left[ e( I J)  - e(I) -e(J)\right]
= e_1(I|J). 
\eeqn
\end{proof}

\section{The adjoint of $I^rJ^s$ and its colength}
In this section, we give  a formula for the colength of the adjoint of $I^rJ^s$, where $I, J$ are $\m$-primary ideals in a two-dimensional   regular local ring. 
This is useful for expressing   the adjoint of $I^rJ^s$ in terms of  the adjoint of $I^r$ and the adjoint of $J^s$. 
We begin by computing the colength of  adjoint of powers of one $\m$-primary ideal. 
\begin{proposition}
\label{prop-hd-one-core}
Let $(R,\m)$ be a two-dimensional regular local ring and $I$ be an $\m$-primary  ideal. Then for all $n \geq 1$, 
\beq
\label{eqn prop-hd-one-core}
   \lambda \left( \frac{R}{\adj{(I^n)}} \right)
= e(I)
    \binom{n+1}{2}
 -     \lambda \left(  \frac{R}{\olin{I}} \right)  
  n.
  \eeq
\end{proposition}
\begin{proof}  
By \cite[Remark~1.2(b)]{lipman-adjoint}, $\adj{(I^n)}$ is a complete ideal  and the point basis of $\adj(I^n)$ is $\max\{ 0, n \cdot r_{_{T}}(I)-1\}_{R \preceq T}$  (\cite[Proposition 3.1.2]{lipman-adjoint}).
   Hence
\begin{align}
\label{hd-oneideal} \nonumber
        &  \lambda \left( \frac{R}{\adj{(I^n)}}  \right)&\\ \nonumber
= &  \sum_{R \preceq T} \binom{n r_{_{T}}(I) }{2}     [T:R]  
  & \mbox{[Theorem~\ref{hoskin-deligne}(\ref{hoskin-deligne-a})]} \\ \nonumber
=&  \left(  \sum_{R \preceq T}{r_{_T}(I)}^2 [T:R] \right) \binom{n+1}{2} 
-       \left(  \sum_{R \preceq T}\binom{r_{_{T}}(I)+1}{2}  
         [T:R] \right)  n &\\
=&   e(I) \binom{n+1}{2}
 -            \lambda\left(\frac{R}{\olin{I}}  \right) n.  & \mbox{[Theorem~\ref{hoskin-deligne}(\ref{hoskin-deligne-a}, \ref{hoskin-deligne-b})]}
\end{align}
\end{proof}
\begin{lemma}
\label{lemma-containment}
Let $R$ be a regular  domain of dimension $d$ and let  $I,K$  be  ideals in  $R$.  Then 
$$K  \cdot \adj(I) \subseteq  \adj(KI).$$
\end{lemma}
\begin{proof}
Let $\upsilon$ be a valuation corresponding to a prime divisor $R_{\upsilon}$ of $R$. Then for any  $a \in K$  and $r \in \adj(I)$, 
\beqn
        \upsilon(ar)
=      \upsilon(a) + \upsilon(r)
\geq \upsilon(K) + \upsilon(I) - \upsilon(J_{R_{\upsilon}/R})
=       \upsilon(KI) - \upsilon(J_{R_{\upsilon}/R}). 
\eeqn
This implies that  $ar\in \adj(KI).$ 
\end{proof}
The next proposition is useful in obtaining our main result. 
\begin{proposition}
\label{hd-two-ideals}
Let $I$ and $J$ be $\m$-primary  ideals in a two-dimensional regular local ring $(R, \m)$.  Then for all $r,s \geq 1$, 
\beq
\label{hd-two-ideals-1}
\lambda\left( \frac{R}{\adj(I^rJ^s)} \right)
&=&  rs~e_1(I|J) 
+  e(I)\binom{r+1}{2} 
-  r \lambda \left( \frac{R}{\olin{I}} \right) 
+  e(J) \binom{s+1}{2} 
-s \lambda \left( \frac{R}{\olin{J}}\right)  \\ \label{hd-two-ideals-2}
&=&  rse_1(I|J)   + \lambda \left( \frac{R}{\adj(I^r)} \right)
+  \lambda \left ( \frac{R}{\adj(J^s)} \right) . 
\eeq
\end{proposition}
\begin{proof} The ideals  $\olin{I^r}$, 
$\olin{J^s}$ and  $ (\olin{I})^r (\olin{J})^s$  are complete ideals (Theorem~\ref{hoskin-deligne}(\ref{hoskin-deligne-0})). 
Moreover,   $e(\olin{I}) = e(I)$,   $e(\olin{J}) = e(J)$  (see \cite{northcott-rees}), and 
$e_1(\olin{I^r} | \olin{J^s}) = rse_{1}( \olin{I}| \olin{J}) = rs e_1(I | J)$  (Theorem~\ref{lipman-length-twoideals}(\ref{lipman-length-twoideals-2},\ref{lipman-length-twoideals-5})).

 By \cite[Remark~1.2(b)]{lipman-adjoint}  $\adj(I^r) = \adj(\ov{I^r})$,  $\adj(J^s) = \adj(\ov{J^s})$  and  $\adj(I^rJ^s) = \adj(\ov{I^rJ^s})$. Hence we can assume that $I$ and $J$ are  complete ideals. Therefore,
\begin{align*}
&~     \lambda\left( \frac{R}{\adj(I^rJ^s)} \right) &\\
=&~ e(I^rJ^s)-\lambda \left( \frac{R}{I^rJ^s} \right)
&\mbox{[Proposition~\ref{prop-hd-one-core}]}\\
=&~ \left[e(I^r) + 2rse_1(I|J)+e(J^s) \right]
- \left[   \lambda \left( \frac{R}{I^r} \right)
+   \lambda \left( \frac{R}{J^s} \right) 
+    rs e_1(I|J) \right] 
&  \mbox{[Lemma~\ref{rees-1961}, Theorem~\ref{lipman-length-twoideals}(\ref{lipman-length-twoideals-1},\ref{lipman-length-twoideals-2})]}\\
=&~ rs e_1(I|J) + e_1(I^r) + e_1(J^s)&\mbox{[Theorem~\ref{hoskin-deligne}(\ref{hoskin-deligne-b})]}\\
=&~  rs e_1(I|J)  + \sum_{R\preceq T} \binom{rr_{_T}(I)}{2} [T:R] + \sum_{R\preceq T} \binom{sr_{_T}(J)}{2} [T:R] 
&\mbox{[Theorem~\ref{hoskin-deligne}(\ref{hoskin-deligne-b})]}\\
=&~ rse_1(I|J) 
+  \left[ e(I)\binom{r+1}{2}-   \lambda \left( \frac{R} {I} \right) r \right]
+ \left[ e(J)\binom{s+1}{2}-  \lambda \left( \frac{R}{J}  \right)  s\right]\\
=&~   rse_1(I|J)  +  \lambda \left( \frac{R}{\adj(I^r)} \right)
+  \lambda \left ( \frac{R}{\adj(J^s)} \right). & \mbox{[Proposition~\ref{prop-hd-one-core}]}
\end{align*}
 \end{proof}
The next result is an analogue of  a result of J.~Lipman and B.~Tessier (see \cite[Corollary~2.2]{verma-nagoya}) and of J. K. Verma  (\cite[Theorem~2.1]{verma-nagoya}) for the adjoint of product of ideals.
\begin{proposition}
\label{joint-reduction-zero}
Let $I$ and $J$ be $\m$-primary ideals in a two-dimensional regular local ring $(R, \m)$.  
\been[\em (a)]
\item
\label{joint-reduction-zero-a}
Let  $\{a,b\}$ is a joint reduction of $\{I, J\}$.   Then for all $r,s\geq 1,$
\begin{align}
\label{joint eqn for adj}
\adj(I^rJ^s)=a^r \cdot \adj(J^s) + b^s \cdot \adj(I^r) = I^r \cdot \adj (J^s) + J^s  \cdot \adj (I^r)  . 
\end{align}

\item
\label{reduction-one-a}
Let $(a,b) $ is a minimal reduction of $I$.   Then for all $n \geq 1$,
$$
\adj(I^{n+1} )=(a,b)  \cdot \adj(I^{n}) = I^{n} \cdot \adj(I).
$$
\eeen
\end{proposition}
\begin{proof}  
(\ref{joint-reduction-zero-a}) 
 By  Lemma~\ref{lemma-containment}  
 \beqn
 a^r \cdot \adj(J^s)
+           b^s \cdot \adj(I^r) 
\subseteq I^r \cdot \adj{J^s}
+         J^s \cdot \adj(I^r) 
\subseteq \adj(I^rJ^s).
\eeqn
To show the  equality, we compute colength of the ideals. 
Since  $\{a,b\}$ is a joint reduction of $\{I,J\}$,   $\{a^r,b^s\}$ is a joint reduction of $\{I^r,J^s\}$ \cite[Theorem~1.3]{Ree84}. Hence from Theorem~\ref{lipman-length-twoideals}(\ref{lipman-length-twoideals-3}) we get
\begin{align}
\label{mix mult of a^r b^s}
e(a^r, b^s) = e_1(I^r |J^s) = rs e_1(I|J). 
\end{align}
 Therefore, 
    \begin{align*}
    \lambda\left( \frac{R}{\adj(I^rJ^s)} \right) 
     = &~e(a^r,b^s)          
     + \lambda\left( \frac{R}{\adj(I^r)}\right) 
     + \lambda\left( \frac{R}{\adj( J^s)}\right)  
     & \mbox{[by (\ref{hd-two-ideals-2})]}\\
    =&~   \lambda \left( \frac{R}{a^r \cdot \adj(J^s) + b^s \cdot \adj(I^r)} \right).  
    & \mbox{\cite[Lemma~3.1]{verma-nagoya}}. 
         \end{align*}

 (\ref{reduction-one-a}) Apply induction on $n$. Since  $(a,b)$ is a minimal reduction of $I$,  $\{a,b\}$ is a joint reduction of $\{I,I\}$.  Put  $J=I$ and $r=s=1$  in (\ref{joint eqn for adj}). We get 
 \begin{align}
\label{red eqn for adj n=1}
\adj(I^2)=a \cdot \adj(I) + b \cdot \adj(I) = I \cdot \adj (I). 
\end{align}
 Hence the result is true for $n=1$. Now let $n \geq 1$. Since $\{a, b^n\}$ is a joint reduction of $\{I, I^n\}$  \cite[Theorem~1.3]{Ree84}. 
 Put $J= I$, $r=1$ and $s = n$ in (\ref{joint eqn for adj}). We get 
 \begin{align*}
 \adj(I^{n+1})
 =&~a \cdot  \adj(I^{n}) + b^{n} \cdot \adj(I)  & \\
 =&~ a (a,b) \cdot \adj(I^{n-1})   + b^n  \cdot \adj(I)  & \mbox{[induction hypothesis]}\\
 \subseteq&~  (a,b)  \adj(I^{n}) & \mbox{[Lemma~\ref{lemma-containment}]}\\
\subseteq&~  \adj(I^{n+1}) .& \mbox{[Lemma~\ref{lemma-containment}]}\\
 \end{align*}
This proves the first equality. 

The second equality is true for $n = 1$ by (\ref{red eqn for adj n=1}). Let $n >1$. 
Put  $J=I$,  $r=1$ and $s=n$. Then from the second equality in  (\ref{joint eqn for adj})  and induction hypothesis we get 
\begin{align*}
       \adj(I^{n+1} )
= I \cdot \adj(I^n) + I^n \cdot \adj(I)
= I  (I^{n-1}  \cdot \adj (I)) + I^n \cdot \adj(I) 
= I^n \cdot \adj(I).  
\end{align*}
\end{proof}

\section{Core of $I^rJ^s$ and its colength}
In this section, we give  a formula for the colength of the core of $I^rJ^s$, where $I, J$ are $\m$-primary  complete ideals in a two-dimensional   regular local ring. 
This is useful for expressing   the core of $I^rJ^s$ in terms of   the core of $I^r$ and the core of $J^s$. 
We begin by computing the colength of  core  of powers of an  $\m$-primary ideal. 
Before we  give  a formula for $\lambda(R/\core(I^n))$, we  prove a few preliminary results.
\begin{lemma}
\label{formula for mult of core}
Let $I$ be a complete ideal  in a two-dimensional regular local ring $(R, \m)$. Then  for all $n,r,s \geq 1$, 
\begin{align}
 \label{mult formla ideal}
          e(I^n) 
=&~ n^2 \left(  \sum_{R \preceq T}  r_{_T}(I)^2 [T:R]\right). \\  \label{mult formla adjoint}
          e(\adj(I^n))
=&~ \sum_{R \preceq T; r_{_T}(I) \geq 1} \left(  (n r_{_T}(I)-1)^2 \right) [T:R]. \\ \label{mult formula ideal and adjoint}
         e(I^r \cdot \adj(I^s))
=&~ \sum_{R \preceq T;  r_{_{T}}(I)  \geq 1}( (r+s) r_{_T}( I) -1 )^2 [T:R].\\ \label{mixmult formula ideal and adjoint}
e_1(I^r|\adj(I^s))
=&~   \sum_{R \preceq T}     r r_{_T}( I)  \left(  s r_{_T}( I)  -1  \right) [T:R].
\end{align}
\end{lemma}
\begin{proof}
By Theorem~\ref{hoskin-deligne}(\ref{hoskin-deligne-b}) we get 
\begin{align*}
           e(I^n) 
= \sum_{R \preceq T} ( r_{_T}(I^n))^2[T:R]
=  \sum_{R \preceq T} (n r_{_T}(I))^2 [T:R] 
= n^2\left (  \sum_{R \preceq T} (r_{_T}(I))^2 [T:R]   \right).
\end{align*}
By \cite[Proposition 3.1.2]{lipman-adjoint}, for all $n \geq 1$, $\adj{(I^n)}$ is a complete ideal  whose point basis is
 $\max\{ 0, r_{_{T}}(I^n)-1\}_{R \preceq T}$.
Hence from Theorem~\ref{hoskin-deligne}(\ref{hoskin-deligne-b}), 
\begin{align*}
           e(\adj(I^n))
=  \sum_{R \preceq T;  r_{_T}(I) \geq 1} \left( r_{_T}(I^n)-1 \right)^2 [T:R]
=      \sum_{R \preceq T; r_{_T}(I) \geq 1} \left(n r_{_T}(I)-1 \right)^2 [T:R].
\end{align*}

Since $I$ and $\adj(I^s)$ are complete ideals, so is $I^r \cdot  \adj(I^s)$ by Theorem~\ref{hoskin-deligne}(\ref{hoskin-deligne-0}). By Proposition~\ref{joint-reduction-zero}(\ref{reduction-one-a}), 
$I^r \cdot  \adj(I^s) = I^r I^{s-1} \cdot \adj(I) = I^{r+s-1} \cdot \adj(I)$. 
By \cite[Proposition 3.1.2]{lipman-adjoint},  the point basis of $I^{r +s-1} \adj(I)$ is 
$\max\{ 0,  r_{_{T}}(I^{r+s-1}) + (r_{_{T}}(I)  -1)\}_{R \preceq T;  r_{_{T}}(I)   \geq 1}$.
\begin{align*}
          e( I^r \cdot \adj(I^s))
 =&~  \sum_{R \preceq T;   r_{_{T}}(I)  \geq 1} \left(  (r+s-1) r_{_{T}}(I) +r_{_{T}}(I)  -1 \right)^2  [T:R] \\
 =&~    \sum_{R \preceq T;  r_{_{T}}(I)  \geq 1}( (r+s) r_{_T}( I) -1 )^2 [T:R].
\end{align*}
For all $r,s \geq 1$,
\begin{align*}
&~       e_1(I^r|\adj(I^s))\\
=&~ \frac{1}{2} \left[  e(I^r \cdot \adj(I^s)) - e(I^r) -e(\adj(I^s))\right] & \mbox{[Lemma~\ref{rees-1961}]}  &\\
=&~ \frac{1}{2} \left[  \sum_{R \preceq T; r_{_T}(I) \geq 1}  \left[   (  (r +  s )r_{_T}( I) -1 )^2 - (r r_{_T}(I))^2 - (s r_{_T}(I)-1)^2 \right] [T:R]  \right]
&      \mbox{[by (\ref{mult formla ideal}),  (\ref{mult formla adjoint}),  (\ref{mult formula ideal and adjoint})]}\\
=&~  \sum_{R \preceq T;  r_{_{T}}(I)  \geq 1}  \left(   rs r_{_T}( I)^2 - r r_{_T}( I)    \right) [T:R] 
& \\
=&~  \sum_{R \preceq T}  \left(   r r_{_T}( I)  \left(  s r_{_T}( I)  -1 \right) \right)[T:R] . 
\end{align*}
\end{proof}

\begin{lemma}
\label{comb lemma}
Let $r,n$ be non-negative integers. Then 
\beqn
r^2 n^2 - nr
= 2 r^2 \binom{n+1}{2} - 2 \binom{r+1}{2} n.
\eeqn
\end{lemma}
\begin{proof} One can verify that 
\beqn
    2 r^2 \binom{n+1}{2} - 2 \binom{r+1}{2} n
&=& \frac{2r^2 n(n+1) - 2 r(r+1) n}{2}\\
&=& r^2 n^2 + r^2n -  r^2n - rn\\
&=& r^2n^2 - rn.
\eeqn
\end{proof}
\begin{proposition}
\label{prop-hd-one}
Let $(R,\m)$ be a two-dimensional regular local ring and $I$ be an $\m$-primary  complete ideal. Then for all $n \geq 1$,
$$
   \lambda \left( \frac{R}{\core(I^n)} \right)
= 4 e(I)
    \binom{n+1}{2}
 -    \left[e(I) +2\lambda\left(\frac{R}{I}  \right) \right]n.
  $$
\end{proposition}
\begin{proof}  
Since $ \core(I^n)= I^n \adj(I^n)$ and  both $I^n$ and $\adj(I^n)$ are complete (Theorem~\ref{hoskin-deligne}(\ref{hoskin-deligne-0})),  $\core(I^n)$ is complete. Hence by  
Theorem~\ref{lipman-length-twoideals}(\ref{lipman-length-twoideals-1}), 
\beq
\label{hd-oneideal-core}
       \lambda \left( \frac{R}{\core(I^n)} \right)
&=& e_1(I^n | \adj(I^n))
+  \lambda \left( \frac{R}{I^n} \right)
+  \lambda \left( \frac{R}{\adj(I^n)} \right).
\eeq
One can verify that 
\begin{align}
\label{hd-oneideal-core-1} \nonumber
&~  e_1(I^n | \adj(I^n)) &  \\ \nonumber
= &~ \sum_{R \preceq T} \left(n^2 r_{_T}(I)^2 -n r_{_T}(I) \right)  [T:R]   & \mbox{[by (\ref{mixmult formula ideal and adjoint})]}   \\ \nonumber
=&~ 2 \left(\sum_{R \preceq T}  r_{_T}(I)^2[T:R] \right)\binom{n+1}{2}
-      2 \left( \sum_{R \preceq T} \binom{r_{_{T}}+1}{2}  
         [T:R] \right)  n & \mbox{[Lemma~\ref{comb lemma}]}\\
 =&~ 2  e(I) \binom{n+1}{2} 
-       2   \lambda\left(\frac{R}{I}  \right) n. 
\end{align}
Substituting (\ref{hd-oneideal-core-1}) in (\ref{hd-oneideal-core})  and applying   Theorem~\ref{hoskin-deligne}(\ref{hoskin-deligne-a},\ref{hoskin-deligne-b})   and 
Proposition~\ref{prop-hd-one-core}
we get 
\beqn
&& \lambda \left( \frac{R}{\core(I^n)} \right)\\
 &=&  2  e(I) \binom{n+1}{2} 
- 2   \lambda\left(\frac{R}{I}  \right) n
+ \left[ e(I) \binom{n+1}{2}  -  \left[e(I) -\lambda\left(\frac{R}{I}  \right) \right]n\right]
+ \left[e(I)    \binom{n+1}{2}-     \lambda \left(  \frac{R}{I} \right)   n\right]\\
  &=& 4  e(I) \binom{n+1}{2} 
  -\left[e(I) +2\lambda\left(\frac{R}{I}  \right) \right]n.
  \eeqn
\end{proof}
In \cite[Proposition~4.4]{hs06}, C.~Huneke and I.~Swanson  proved that $I^{2n-1}  \cdot \adj(I) = \core(I^{n})$ for all $n \geq 1$. We recover this result. 
\begin{corollary}
\label{core-huneke-swanson}
    Let $(R,\m)$ be a two-dimensional regular local ring and $I$ be $\m$-primary complete ideal in a $R.$ Then for all $n \geq 1$, 
    \been[\em (a)]
    \item
    \label{core-huneke-swanson-2}
    $ \core(I^n)=I^{2n-1} \cdot  \adj(I) = \adj(I^{2n}). $
    \item
    \label{core-huneke-swanson-3}
      $\core(I^{n+1}) =I^{2} \cdot \core (I^n).$
      \eeen
\end{corollary}
    \begin{proof} (\ref{core-huneke-swanson-2})  
  Applying Lemma~\ref{lemma-containment} we get,   $\core(I^n) = I^n \cdot \adj(I^n) \subseteq \adj(I^{2n})$  and $I^{2n-1} \cdot \adj(I) \subseteq I^n \cdot \adj(I^n) = \core(I^n)$. 
  
  Hence to show  the equality we  compute the respective colengths. 
 Since  both $I^{2n-1}$ and $\adj(I)$ are complete, by   Theorem~\ref{lipman-length-twoideals}(\ref{lipman-length-twoideals-1}) we get 
\beq
\label{hd-oneideal-core-2n-1adj(I)}
       \lambda \left( \frac{R}{I^{2n-1}   \cdot \adj(I)} \right)
&=& e_1(I^{2n-1} | \adj(I))
+  \lambda \left( \frac{R}{I^{2n-1}} \right)
+  \lambda \left( \frac{R}{\adj(I)} \right).
\eeq
One can verify that 
\begin{align}
\label{hd-oneideal-core-2} \nonumber
       e_1(I^{2n-1} | \adj(I))
= &~ \sum_{R \preceq T} \left((2n-1)r_{_T}(I) (r_{_T}(I)-1)\right) [T:R]    &\mbox{ [by (\ref{mixmult formula ideal and adjoint})]} \\ \nonumber
=&~ 2 ( 2n-1) \sum_{R \preceq T} \binom{r_{_T}(I)}{2} [T:R]\\
=&~ 2(2n-1) \left[ e(I) - \lambda\left(\frac{R}{I} \right) \right].
& \mbox{[Theorem~\ref{hoskin-deligne}(\ref{hoskin-deligne-b})]}
\end{align}
Substituting (\ref{hd-oneideal-core-2}) in   (\ref{hd-oneideal-core-2n-1adj(I)})   and applying  Theorem~\ref{hoskin-deligne}(\ref{hoskin-deligne-a},\ref{hoskin-deligne-b}) and (\ref{eqn prop-hd-one-core}) we get 
\beq
\label{hd for core-1} \nonumber
 &&      \lambda \left( \frac{R}{I^{2n-1}  \adj(I) } \right)\\ \nonumber
&=& 2(2n-1) \left[ e(I) -\lambda\left(\frac{R}{I} \right) \right]
+  \left[e(I) \binom{2n}{2}  -  \left[e(I) -\lambda\left(\frac{R}{I}  \right) \right](2n-1)\right]
+ \left[e(I)-     \lambda \left(  \frac{R}{I} \right) \right] \\
&=& 4  e(I) \binom{n+1}{2} 
  -\left[e(I) +2\lambda\left(\frac{R}{I}  \right) \right]n. 
  \eeq
From (\ref{eqn prop-hd-one-core}) we get 
\begin{align}
\label{length of adj(I^{2n})} \nonumber
   \lambda \left( \frac{R}{\adj(I^{2n})} \right)
= &~e(I) \binom{2n+1}{2}
 -     \lambda \left(  \frac{R}{I} \right)  (2n)\\ \nonumber
 = &~ 4e(I) \binom{n+1}{2} - e(I)n -  2\lambda \left(  \frac{R}{I} \right) n\\
 = &~ 4e(I) \binom{n+1}{2} - \left[e(I) -  2\lambda \left(  \frac{R}{I} \right) \right]n.
  \end{align}
From (\ref{hd for core-1}),  (\ref{length of adj(I^{2n})})  and Proposition~\ref{prop-hd-one} the result follows. 

(\ref{core-huneke-swanson-3})
From (\ref{core-huneke-swanson-2}) we get that for all $n \geq 1$,
\beq
 \core(I^{n+1}) 
=  I^{2(n+1)-1}  \cdot \adj{I}
= I^2 (I^{2n-1}  \cdot \adj{I})
=I^{2} \core (I^n).
\eeq
    \end{proof}
We now calculate the colength of $\core(I^rJ^s)$ using the Hoskin-Deligne formula.
 \begin{proposition}
\label{hd-two-ideals core prop}
Let $I$ and $J$ be $\m$-primary complete ideals in a two-dimensional regular local ring $(R, \m)$.  Then for all $r,s \geq 1$, 
\beqn
\lambda\left( \frac{R}{\core(I^rJ^s)} \right) 
= 4 rs e_1(I|J)
+  \lambda\left( \frac{R}{\core(I^r)} \right)
+  \lambda\left( \frac{R}{\core(J^s)} \right). 
\eeqn
\end{proposition} 
\begin{proof} Using Proposition~\ref{prop-hd-one} we get that for all $r,s \geq 1$, 
\begin{align*}
 &~  \lambda \left( \frac{R}{\core(I^rJ^s)} \right)\\
=&~ 4 e(I^rJ^s)
    -    \left[e(I^rJ^s) +2\lambda\left(\frac{R}{I^rJ^s}  \right) \right]\\
=&~ 3 e(I^rJ^s)-  2  \lambda\left(\frac{R}{I^rJ^s}  \right) \\ 
 =&~ 3 \left[e(I^r) +2 e_1(I^r|J^s) + e(J^s)\right]
 -      2 \left[  e_1(I^r | J^s) 
 +        \lambda\left(\frac{R}{I^r}  \right) 
 +        \lambda\left(\frac{R}{J^s}  \right) \right]    
  &      \mbox {[Lemma~\ref{rees-1961}, Theorem~\ref{lipman-length-twoideals}(\ref{lipman-length-twoideals-1})]}\\
 =&~ 4 rs e_1(I| J)   + 3 e(I^r) + 3 e(J^s) - 2   \lambda\left(\frac{R}{I^r}  \right) -2  \lambda\left(\frac{R}{J^s}  \right)\\
 =&~ 4 rs e_1(I| J)   + 3 e(I^r) + 3 e(J^s)  & \\
  &~ +  \left[ -3 e(I^r) +  \lambda \left( \frac{R}{\core(I^r)} \right)  \right]
  +   \left[ -3 e(J^s) +  \lambda \left( \frac{R}{\core(J^s)} \right)  \right] & \mbox{[Proposition~\ref{prop-hd-one}]} \\
 =&~ 4 rs e_1(I| J) 
 +  \lambda \left( \frac{R}{\core(I^r)} \right)
 +  \lambda \left( \frac{R}{\core(J^s)} \right).  
  \end{align*}
\end{proof}
\begin{lemma}
\label{lemma-containment-core}
Let $R$ be a regular  domain of dimension $d$ and let  $I,J$  be complete ideals in $R$.  Then 
$J^2 \core(I) \subseteq  \core(JI)$. 
\end{lemma}
\begin{proof}
Applying Lemma~\ref{lemma-containment} we get 
\beqn
J^2 \cdot \core(I)
=J  I ( J \cdot \adj(I))
\subseteq J  I (\adj{JI})
= \core(JI).
\eeqn
\end{proof}

\begin{theorem}
\label{joint-reduction-zero-core}
Let $I$ and $J$ be $\m$-primary ideals in a two-dimensional  regular local ring $(R, \m)$.   
\been[\em (a)]
\item
\label{joint-reduction-zero-a-core}
Let $\{a,b\}$ be a joint reduction of $\{I,J \}$. Then for all $r,s\geq 1,$
 $$
 \core(I^rJ^s) = a^{2r} \core(J^s) + b^{2s} \core(I^r) = I^{2r} \core(J^s) + \core(I^r) J^{2s}.
 $$
 
\item
\label{reduction-one-a-core}
Let $(a,b)$ be a minimal reduction of $I$. Then all $n \geq 1$, 
$\core(I^{{n+1}}) = (a^2,b^2) \core (I^{n})$. 
\eeen
\end{theorem}
\begin{proof}
(\ref{joint-reduction-zero-a-core})    
From Lemma~\ref{lemma-containment-core}, 
$               a^{2r} \cdot \core(J^s) 
\subseteq I^{2r} \cdot \core(J^s) 
\subseteq \core( I^rJ^s )$. Similarly, 
     $b^{2s} \core(I^r) \subseteq \core(I^rJ^s).$ 
     Therefore, $a^{2r} \core(J^s) + b^{2s} \core(I^r) \subseteq  I^{2r} \core(J^s) + J^{2s} \core(I^r) \subseteq \core(I^rJ^s).$
 To show the equality we compute colength.

Since $\{a,b\}$ is a joint-reduction of $\{I,J\}$, $\{a^{2r}, b^{2s}\}$ is a joint reduction of $\{I^{2r}, J^{2s}\}$. Hence by Theorem~\ref{lipman-length-twoideals}(\ref{lipman-length-twoideals-2}), (\ref{lipman-length-twoideals-3}), (\ref{lipman-length-twoideals-4}) we get 
\beq 
\label{mixed mult of powers}
4rse_1(I|J)  = e_1(I^{2r}| J^{2s}) = e(a^{2r}, b^{2s}). 
\eeq

Hence from   Proposition~\ref{hd-two-ideals core prop} we get 
  \begin{align}  \nonumber
   \label{eqn-1-joint-reduction-zero-a-core}
         \lambda\left( \frac{R}{\core(I^rJ^s)} \right) 
 =&~ 4rse_1(I| J)          
 +     \lambda\left( \frac{R}{\core(I^r)}\right) 
 +    \lambda\left( \frac{R}{\core(J^s)}\right) &  \\  \nonumber
= &~  e(a^{2r},b^{2s})          
 +      \lambda\left( \frac{R}{\core(I^r)}\right) 
 +     \lambda\left( \frac{R}{\core(J^s)}\right) &  \\  \nonumber
=&~       \lambda \left( \frac{R}
                                     { a^{2r} \cdot \core(J^s) + b^{2s}  \cdot \core(I^r) } \right)  
&    \mbox{\cite[Lemma~3.1]{verma-nagoya}} \\ \nonumber
\geq &~ \lambda\left( \frac{R}{  I^r \cdot \core(J^s) +J^s \cdot \core(I^r)   } \right). &\\
\geq &~ \lambda\left( \frac{R}{ \core(I^rJ^s) } \right). &
 \end{align}
  Hence equality holds in (\ref{eqn-1-joint-reduction-zero-a-core}) which proves (\ref{joint-reduction-zero-a-core}). 
 
 (\ref{reduction-one-a-core}) 
 Since $(a,b)$ is a reduction of $I$, $\{a, b^n\}$ is a joint reduction of $\{I, I^n\}$. Hence by (\ref{joint-reduction-zero-a-core})  we have
 \begin{align*}
        \core(I^{n+1})
  =&~ \core( I I^n)\\
  =&~ a^2 \core(I^n) + b^{2n} \core(I) &  \mbox{[by (\ref{joint-reduction-zero-a-core})]}\\
  \subseteq &~  a^2 \core(I^n) + b^2 I^{2(n-1)}\core(I)& \\
  \subseteq&~ a^2 \core(I^n) + b^2 \core (I^n) & \mbox{[Lemma~\ref{lemma-containment-core}]}\\
  \subseteq &~ \core(I^{n+1}). & \mbox{[Lemma~\ref{lemma-containment-core}]}
 \end{align*}
\end{proof}

 \section{Examples}
\begin{example} 
{\rm
Let $R =\mathbb{Q}[[x,y]]$ be the formal power series ring over $\mathbb{Q}$ and $\m = (x,y)$. Let   $I=(x^2,xy,y^3)$  and $ K=(x^3,xy,y^2)$. For all $r,s \geq 1$, we  compute  
$\adj(I^r K^s)$,  $\core(I^rK^s)$,  and their respective colength.

The ideals  $(x,y^3)$  and  $(y,x^2)$ are complete ideals. Therefore $I=(x,y^3)\cap (y,x^2)$  is a complete ideal.  Similarly, $K$ is a complete ideal. 
Put 
$ J=(xy,x^2+y^3)$ and $L=(xy,y^2+x^3).$
Since 
 $JI=I^2$ (resp.  $LK=K^2$), $J$  (resp. $L$) is a reduction of $I$ (resp. $K$). 
By Lech's formula \cite{lech}
\begin{align*}
    e(xy,x^2+y^3)
=&~ e(x,x^2+y^3)+e(y,x^2+y^3)
    = e(x, y^3) + e(y, x^2)
    = 3+2 
    =5.\\
        e(xy,x^3+y^2)
    =&~5.  & \mbox{[by symmetry]}
\end{align*}
Hence 
\beq
\label{formula for multiplicity}
e(I)=e(J)=e(K)=e(L)=5.
\eeq
Since,   $  IK=y^2I+x^2K,$
$ \{x^2,y^2\}$ is a joint reduction of $\{I,K\}$ we get  
\beq
\label{joint red of I and K}
e_1(I|K)=e(x^2, y^2)= 4. \hspace{.2in} \mbox{[Theorem~\ref{lipman-length-twoideals}(\ref{lipman-length-twoideals-3})]}
\eeq
Therefore  for all $r,s \ge 1$
\begin{align*}
    e(I^rK^s)  = e(I)r^2+2 e_1(I|K)rs+e(J)s^2  = 5r^2+8rs+5s^2. \hspace{.2in}  \mbox{[Lemma~\ref{rees-1961}, Theorem~\ref{lipman-length-twoideals}(\ref{lipman-length-twoideals-2})]}
\end{align*}

(a)   We claim  that for all $n \geq 1$
\begin{align}
\label{adjoint of In}
 \adj(I^n)=\m I^{n-1} & \mbox{ and }  \adj(K^n)=\m K^{n-1}.
 \end{align}
We  prove  the claim only for the ideal $I$ as the proof for $K$ is similar. 
Apply induction on $n.$ Let  $n=1.$ As $J$ is a minimal reduction of $I$, applying \cite[Proposition 3.3]{lipman-adjoint} we get 
$\m \subseteq J:I = \adj(I)$. To show $\adj(I) \subseteq  \m$, we consider the point basis of $I$. Note that $r_{R}(I) = 2$ and for $T = R[\m/y]_{ (x/y,y)}$, $r_T(I) =1$. Hence,  the point basis of $\adj(I)$ is:  $\{ r_R(I)=1,0\}$ which implies that $\adj(I) \not =R$ and therefore 
$\adj(I) = \m$.
Hence the claim  is true for $n=1$.

Now let $n>1$.  Applying   \cite[(2.3)]{lipman-adjoint}  and induction hypothesis we get 
\begin{align*}
    \adj(I^{n}) = I \cdot \adj(I^{n-1}) = I (\m I^{n-2})=\m I^{n-1},
\end{align*}
which proves (a). 

(b) We compute $\lambda(R/\adj(I^n)).$
 As $\m$ and $I$ are complete ideals, the ideals $I^{n-1}$ and $\m I^{n-1}$    are complete (Theorem~\ref{hoskin-deligne}(\ref{hoskin-deligne-0})). Let   
$\mu(I^{n-1}):=\lambda(I^{n-1}/\m I^{n-1})$ denote the minimal number of generators of $I^{n-1}$. Then
\begin{align}
\label{number of generators} 
        \mu(I^{n-1})
=&~   1+r_{_R} ( I^{n-1})=  1 + 2(n-1) = 2n-1,  & \mbox{\cite[Proposition~2.3]{huneke}}\\ \label{number of generators-1} 
          e_1(I) 
 = &~ e(I) - \lambda\left( \frac{R}{I}\right)  = 5 - 4 = 1. 
& \mbox { [Theorem~\ref{hoskin-deligne}(\ref{hoskin-deligne-b}), (\ref{formula for multiplicity})]}
\end{align}
 Hence, for all $n\ge 1,$ 
\begin{align}
\label{length of adj(I)} \nonumber
    \lambda \left(  \frac{R}{\adj(I^n) } \right)
=&~ \lambda \left( \frac{R}{\m I^{n-1}} \right) & \mbox{[by (\ref{adjoint of In})]}\\ \nonumber
=&~ \lambda \left( \frac{R}{I^{n-1}} \right)
 +    \lambda \left(\frac{I^{n-1}}{\m I^{n-1}}\right) \\ \nonumber
 = &~  5\binom{n}{2}-(n-1)+2n-1 
 &     \mbox{[Theorem~\ref{hoskin-deligne}(\ref{hoskin-deligne-a}), (\ref{formula for multiplicity}),  (\ref{number of generators}) (\ref{number of generators-1})]}\\
=&~  5 \binom{n+1}{2} -4n.
\end{align}

\noindent (c) We  describe  $\adj(I^rK^s)$ and compute $\lambda(R/ (\adj(I^r K^s))$. 
For all   $r,s \ge 1$ 
\begin{align}
\label{formula adjoint products} \nonumber
       \adj(I^rK^s) 
= &~  I^r \cdot \adj(K^s) + K^s \cdot \adj(I^r)   & \mbox{[by (\ref{joint eqn for adj})]}\\ \nonumber
= &~ I^r \m K^{s-1} + K^s \m I^{r-1}  & \mbox{[by (\ref{adjoint of In})]}\\ \nonumber
= &~ \m I^{r-1}K^{s-1} (I+K) \\
= &~ \m^3 I^{r-1}K^{s-1}.
\end{align}
As $r_R(\m)=1$,  $r_R(I) = r_R(K) = 2$ and $r_T(\m) =0$ for all $R \prec T $, 
\begin{align}
\label{mix mult of m and I-1}
      e_1(\m | I) 
=&~ e_1(\m |K)=2  
&   \mbox{[Theorem~\ref{hoskin-deligne}(\ref{hoskin-deligne-c})]}\\ \label{mix mult of m and I-2}
        e_1( \m^3 | I^{r-1}) 
= &~ 3(r-1) e_1(\m|I) = 6(r-1) 
& \mbox{[Theorem~\ref{lipman-length-twoideals}(\ref{lipman-length-twoideals-2}), (\ref{mix mult of m and I-1})]}\\  \label{mix mult of m and I-3}
       e_1( \m^3 | K^{s-1}) 
 =&~  3(s-1) e_1(\m|K) = 6(s-1) 
 & \mbox{[Theorem~\ref{lipman-length-twoideals}(\ref{lipman-length-twoideals-2}), (\ref{mix mult of m and I-1})]}\\  \label{mix mult of m and I-4}
e_1( I^{r-1} | K^{s-1}) 
= &~ (r-1) (s-1) e_1(I| K) = 4 (r-1)(s-1). &  \mbox{[Theorem~\ref{lipman-length-twoideals}(\ref{lipman-length-twoideals-2}), (\ref{joint red of I and K})]}
\end{align}
For all $r,s \geq 1$, 
\begin{align}
\label{length formula for IrKs} \nonumber
   &~ \lambda \left( \frac{R}{ \m^3 I^{r-1} K^{s-1}} \right)  &  \\ \nonumber
  =&~ \lambda \left( \frac{R}{ \m^3 I^{r-1}} \right) 
 +       \lambda \left( \frac{R}{K^{s-1}}\right)    
 +      e_1(\m^3 I^{r-1}|K^{s-1}) 
 & \mbox{[Theorem~\ref{lipman-length-twoideals}(\ref{lipman-length-twoideals-1})]}  & \\ \nonumber
=&  \left[\lambda \left( \frac{R}{\m^3}\right) 
 + e_1(\m^3|I^{r-1}) 
 + \lambda\left( \frac{R}{I^{r-1}}\right)  \right] &  \\  \nonumber
& + \lambda \left (\frac{R}{K^{s-1}}\right)
   + \left[e_1(\m^3|K^{s-1}) + e_1(I^{r-1}|K^{s-1}) \right]   
 &     \mbox{[Theorem~\ref{lipman-length-twoideals}(\ref{lipman-length-twoideals-1}, \ref{lipman-length-twoideals-4}))]} \\ \nonumber 
  =&~ 6+ 6(r-1) +  \left[5  \binom{r}{2}  - (r-1)\right]  + \left[5\binom{s}{2}  - (s-1)\right] \\ \nonumber
&     + 6(s-1) +4(r-1)(s-1)   
&   \mbox{[Theorem~\ref{hoskin-deligne}(\ref{hoskin-deligne-a}),  (\ref{mix mult of m and I-2}), (\ref{mix mult of m and I-3}), (\ref{mix mult of m and I-4})]} 
        \\ \nonumber
 =&~  6r + \left[5 \binom{r+1}{2}  - 5r-(r-1)\right]    + \left[5 \binom{s+1}{2}  - 5s-(s-1)\right]\\ \nonumber
    & + 6s-6 + 4rs -4r-4s + 4 \\  
=&~ 5 \binom{r+1}{2}  + 4rs + 5 \binom{s+1}{2}  -4r -4s. 
\end{align} 

\noindent (d) We compute $\core (I^rK^s)$ and
$\lambda (R/\core(I^rJ^s)).$
Applying (\ref{formula adjoint products}) we get
\begin{align*}
    \core (I^rK^s) & = I^rK^s \cdot \adj(I^rK^s )  = I^rK^s \m^3 I^{r-1}K^{s-1}  = \m^3 I^{2r-1}K^{2s-1}.
\end{align*}
Replacing $r-1$ by $2r-1$ and $s-1$ by $2s-1$ in (\ref{length formula for IrKs}) we get  
\begin{align*}
   & \lambda \left( \frac{R}{ \m^3 I^{2r-1} K^{2s-1}} \right) \\ \nonumber
    =& 5 \binom{2r+1}{2}  + 16rs + 5 \binom{2s+1}{2}  -8r -8s\\ \nonumber
    =&  5 \left[ 4 \binom{r+1}{2} -r \right] + 16rs + 5\left[4 \binom{s+1}{2}-s\right]  -8r -8s\\
= & 20 \binom{r+1}{2} + 16rs +20 \binom{s+1}{2} -13r-13s.
\end{align*} 
} 
\end{example}
We extend the previous example and use alternative methods to the compute  the adjoint and the core. 
We will use the following algorithm which  was  given in the thesis of M.~J.~Rhodes for three generated monomials \cite[Lemma~6.5]{Rhodes} and was proved in a more general setting for ideals in a regular domain
\cite[Proposition 18.3.2]{huneke-swanson}. 
\begin{principle}
Let $R$ be a regular domain. Let an ideal $\mathfrak{m}$ of height $h$ be minimally generated by $(x_1, \ldots, x_h)$ and assume that $R/(x_1, \ldots, x_h)$ is regular. Then for any ideal $I$ in $R$,
\beq
\label{rhodes formula}
\adj( I) 
= \bigcap_{i=1}^{h} 
 \left(   \frac{1}{x^{h-1}_i} 
    \adj \left( I R\left[ \frac{\m}{x_i} \right] \right) \cap R \right). 
\eeq
\end{principle}
We recall the following remark of J.~Lipman: 
\begin{remark}
\label{lipman 1.2}
\cite[Remark~1.2(c)]{lipman-adjoint}
Let $R$ be a regular normal domain and $I$ and ideal in $R$. For any $r \in R$, 
$\adj(r I) = r \adj(I)$. 
\end{remark}
We prove a lemma which will be  used in the example. 
\begin{lemma}
Let  $(R, m)$  be a  two-dimensional regular local ring with $\m = (x, y)$.  Let $n \geq 2$ and 
put 
\beq
\label{defn o f Kn}
K_{n}(x,y) :=   (x^{2n-3}, x^{n-2}y, y^{2}  ).
\eeq
Then
  \beq
  \label{adj pf Kn}
\adj ( K_{n}(x,y) ) =   (x^{n-2},y).
\eeq
\end{lemma}
\begin{proof} Apply induction on $n$. 
If $n=2$, then 
$
\adj (K_{2}(x,y)) = \adj( (x,y) ) = R.
$
Now let $n>2$. Then 
\begin{align*}
  &          \adj(K_n (x,y)) \\ 
= &~      \left(\frac{1}{x} \cdot \adj (K_{n}(x,y)) R\left[\frac{\m}{x}\right] \cap R\right)  
\bigcap  \left(\frac{1}{y} \cdot \adj (K_{n}(x,y)) R\left[\frac{\m}{y}\right] \cap R\right) 
&           \mbox{[by (\ref{rhodes formula})]}       \\ 
=&~       \left(  \frac{1}{x} \cdot \adj \left(  x^{2n-3}, x^{n-1} \frac{y}{x}, x^2  \left(\frac{y}{x}\right)^2 \right)    R \left[ \frac{\m }{x}\right]
               \cap R \right)\\
&           \bigcap \left(  \frac{1}{y} \cdot \adj \left(  \left( \frac{x}{y}\right)^{2n-3} y^{2n-3},  
                           \left( \frac{x}{y}\right)^{n-2} y^{n-1},   y^2\right)     R \left[ \frac{\m }{y}\right] \cap R \right) 
&               \\ 
=&~      \left(  \frac{x^2}{x} \cdot \adj \left(K_{n-1} ( x, y/x)  \right)R \left[ \frac{\m }{x} \right] \cap R \right)
\bigcap \left( \frac{y^2}{y}  R \left[ \frac{\m }{y} \right] \cap R\right) 
&           \mbox{[Remark~\ref{lipman 1.2}]}\\ 
=&~       \left( x \left(x^{n-3}, \frac{y}{x}\right) R \left[ \frac{\m }{x} \right] \cap R\right) \cap \m \hspace{.2in} &  \mbox{[by induction hypothesis]}\\ 
=&~      \left( (x^{n-2},  y) R \left[ \frac{\m}{x} \right] \cap R \right) \cap \m & \\ 
=&~      (x^{n-2},y) \cap \m & \\
=&~        (x^{n-2},y) . 
\end{align*}
\end{proof}
\begin{remark}
Since Example~\ref{example general} involves monomial ideals, we can compute the adjoint using the method given in \cite{hubl-swanson}. 
\end{remark}

\begin{example}
\label{example general}
{\rm
Let $R = \mathbb{Q}[[x,y]]$ be the power series ring over $\mathbb{Q}$ and $\mathfrak{m} = (x,y)$. 
Let $u \geq 3$  be an integer. 
Let  $I = (x^{u}, xy, y^{u+1})$ and $K = (x^{u+1}, xy, y^{u})$. For all $r,s\geq 1$, we  compute  
$\adj(I^r K^s)$,  $\core(I^rK^s)$,  and their respective colength. 

Since   $I = (x,y^{u+1}) \cap (y,x^{u})$ and  the ideals  $(x,y^{u+1})$  and  $(y,x^{u})$ are complete ideals. Therefore, 
$I$ is a complete ideal. Similarly, $K$ is a complete ideal. 

Put $J = (x^{u} + y^{u+1}, xy)$ and  $L = (x^{u+1} + y^{u}, xy).$
We  show that $J$ is a reduction of $I$. 
It is enough to show that $I^2 \subseteq JI$ as the other inclusion is obvious. 
Since $(xy)^2 \in JI$,  $xy^{u+2} = (xy) (y^{u+1}) \in JI \subseteq I^2$ and as $u \geq 3$,  
$x^{u} y^{u+1} = (x^2y^2)(x^{u-2}y^{u-1}) \in JI$.
Hence,  
\begin{align*}
     I^2&=(x^{2u},(xy)^2,y^{2(u+1)},x^{u+1}y,xy^{u+2},x^{u} y^{u+1}) = (x^{2u},(xy)^2,y^{2(u+1)},x^{u+1}y,xy^{u+2})\\
    JI &= (x^{u} + y^{u+1}, xy)(x^{u}, xy, y^{u+1}) \\
       &= (x^{2u}+x^{u} y^{u+1}, x^{u+1}y+xy^{u+2},x^{u} y^{u+1}+y^{2(u+1)},x^{u+1}y,(xy)^2,xy^{u+2})\\
       &=  (x^{2u}, x^{u+1}y, y^{2(u+1)},  (xy)^2, xy^{u+2}), 
\end{align*}
This implies that  $ I^2 = JI$.   A similar computation will show that  $L$ is a  reduction of $K$.
By Lech's formula  \cite{lech}
\begin{align*} 
    e(xy, x^{u}+y^{u+1})
    =& e(x,x^{u}+y^{u+1})+e(y,x^{u}+y^{u+1})
    = u+1+u 
    =2u+1 &\\
        e(xy,x^{u+1}+y^{u})
   =&  2u+1.& \mbox{[by symmetry]}
\end{align*}
Hence 
\beq
\label{multiplicity formula general}
e(I)=e(J)=e(K)=e(L)=2u+1.
\eeq
The exact sequence 
\beqn
0 \rightarrow \frac{R}{I} 
\rightarrow \frac{R}{(x,y^{u+1})} \oplus  \frac{R}{(x^{u},y)}
\rightarrow   \frac{R}{(x,y)}
\rightarrow 0
\eeqn
implies that 
\beq
\label{length for r general}
   \lambda \left( \frac{R}{I}  \right)
= \lambda \left(  \frac{R}{(x, y^{u+1})} \right) + \lambda \left(   \frac{R}{(x^{u},y)}\right)
-   \lambda \left( \frac{R}{(x,y)}\right)
= u+1 + u -1 = 2u
\eeq
We claim  that 
\begin{align}
\label{joint reduction general}
IK  = (xy + x^{u}) K + (xy + y^{u}) I.
\end{align}
  One can verify that 
\beq
\label{eqn for IK general}
IK = ( x^{2u+1}, x^{u+1}y,   x^2y^2, xy^{u+1} , y^{2u+1}).
\eeq
As $1 -x^{u-2}y^{u-2}$ is a unit, from the equation 
\begin{align*}
        x^2 y^2 (  1 -x^{u-2} y^{u-2}) 
= - (xy + y^{u}) x^{u}  +  (xy + x^{u}) xy,  
\end{align*}
we get $x^2 y^2 \in  (xy + x^{u}) K + (xy + y^{u}) I.$
Hence, from the equations 
\begin{align*}
x^{2u+1}
&=   -x(xy+ x^{u}) xy 
   + (xy + x^{u})x^{u+1}  + x( x^2y^2)\\ 
x^{u+1}y 
&=    (xy + x^{u}) xy  -x^2y^2 \\
  xy^{u+1}
&=  (xy + y^{u}) xy - x^2 y^2\\ 
    y^{2u+1}
&=  -y(xy+ y^{u}) xy + (xy + y^{u})y^{u+1} + y(x^2y^2)
 \end{align*}
we get $IK \subseteq (xy + x^{u}) K + (xy + y^{u}) I$. The other inclusion is obvious. 
Hence 
$\{    (xy + x^{u}), (xy + y^{u}) \}$ is a joint reduction of $\{I, K\}$. 
Since $1 -x^{u-2}y^{u-2}$ is a unit and
\begin{align*}
    x (1-x^{u-2} y^{u-2} )  =  & x+ y^{u-1} - y^{u-2} (y + x^{u-1}),
\end{align*}
we get 
\begin{align}
\label{mult of y+ x^{u-1}, x+ y^{u-1}}
 e(y+ x^{u-1}, x+ y^{u-1}) = e(x,y)=1
 \end{align}
By using Lech's formula \cite{lech} and   (\ref{mult of y+ x^{u-1}, x+ y^{u-1}}) we get 
\begin{align}
\label{joint red of I and K general} \nonumber
  e_1(I|K)
&= e(xy + x^{u}, xy + y^{u})\\ \nonumber
&= e(x, y) + e(x, x+y^{u-1}) + e( y+ x^{u-1},y) + e(y+ x^{u-1}, x+ y^{u-1})\\ \nonumber
&= 1 + u-1 + u-1 + 1\\
&= 2u.  
\end{align}

\noindent
(a) We claim that $\adj(I^n)=\m I^{n-1}$ and $\adj(K^n)=\m K^{n-1}$   for all $n \ge 1.$
We  prove  the claim only for the ideal $I$ as the proof for $K$ is similar. 
Put $S_1 = R[\m / x]$ and $S_2 = R[\m /y]$. Then applying (\ref{rhodes formula}) we get
\begin{align*}
\adj(I)
=&~    \left(  \left( \frac{1}{x} \cdot \adj(I S_1)  \right)\cap R\right) 
\cap    \left( \left(\frac{1}{y} \cdot \adj(I S_2) \right)\cap R\right) \\
=&~ \left( \frac{1}{x} \adj  \left( (x^{u}, x^2 \frac{y}{x}, x^{u+1} \left(\frac{y}{x} \right)^{u+1} ) S_1 \right) \cap R \right)\\
&\cap  \left( \frac{1}{y}  \adj \left( ( y^{u} \left(  \frac{x}{y} \right)^{u},  \left( \frac{x}{y} \right) y^2, y^{u+1} ) S_2 \right) \cap R \right)\\
=&~   \left(x \cdot \adj \left( x^{u-2}, \frac{y}{x}  \right) S_1  \cap R\right)
\cap  \left( y \cdot \adj \left(  \frac{x}{y} , y^{u-2}  \right) S_2 \cap R\right)  & \mbox{[Remark~\ref{lipman 1.2}]}\\
=&~ (x S_1 \cap R) \cap (y S_2 \cap R)\\
=&~ \m.
\end{align*}
Now let $n>1$. Then applying   \cite[(2.3)]{lipman-adjoint}  and by induction hypothesis we get 
\begin{align}
\label{formula for adj general}
    \adj(I^{n}) = I \cdot \adj(I^{n-1}) = I (\m I^{n-2})=\m I^{n-1},
\end{align}
which proves (a). 

\noindent
(b) We compute $\lambda(R/\adj(I^n)).$
As $\m$ and $I$ are complete ideals, the ideals $I^{n-1}$ and $\m I^{n-1}$    are also  complete (Theorem~\ref{hoskin-deligne}(\ref{hoskin-deligne-0})). Let   
$\mu(I^{n-1}):=\lambda(I^{n-1}/\m I^{n-1})$ denote the minimal number of generators of $I^{n-1}$. Then
\begin{align}
\label{number of generators general} 
       \mu(I^{n-1})
=&~   1+r_{_R} ( I^{n-1})=  1 + 2(n-1) = 2n-1,  & \mbox{\cite[Proposition~2.3]{huneke}}\\ \label{first coefficient r}
     e_1(I) 
= &~ e(I) - \lambda\left( \frac{R}{I}\right)  = (2u+1) - 2u = 1. 
& \mbox { [Theorem~\ref{hoskin-deligne}(\ref{hoskin-deligne-b}), (\ref{multiplicity  formula general}), (\ref{length for r general})]}
\end{align}
Hence, for all $n\ge 1,$ 
 \begin{align}
\label{length of adj(I) general} \nonumber
         \lambda \left(  \frac{R}{\adj(I^n) } \right)
 =&~ \lambda \left( \frac{R}{\m I^{n-1}} \right) & \mbox{ [by  (\ref{formula for adj general})]}\\ \nonumber
 =&~ \lambda \left( \frac{R}{I^{n-1}} \right)
 +  \lambda \left(\frac{I^{n-1}}{\m I^{n-1}}\right) \\ \nonumber
 = &~  (2u+1)\binom{n}{2}-(n-1)+2n-1
 & \mbox{[Theorem~\ref{hoskin-deligne}(\ref{hoskin-deligne-a}), (\ref{multiplicity  formula general}),  (\ref{number of generators general}),  (\ref{first coefficient r})]}\\ \nonumber
     =& \left[(2u+1) \binom{n+1}{2} - (2u+1)n\right] +  n\\ 
      =& (2u+1) \binom{n+1}{2} - (2u)n.
\end{align}

\noindent (c) We compute $\adj(I^rK^s) $. 
We first compute $\adj(IK)$. 
Since  $u \geq 3$, $1+x^{u-3}y^{2u-3}$ and $1-x^{u-2}y^{u-2}$ are units and  the equations 
\beqn
x^2y (1+x^{u-3}y^{2u-3}) 
&=& (1+x^{u-2}y^{u-2}) (xy+y^{u})x -    y^{u-2}(xy+x^{u})y\\
xy^2(1-x^{u-2}y^{u-2}) 
&=& (xy+x^{u})y-x^{u-2}(xy+y^{u})x,
\eeqn
imply that $x^2y, xy^2 \in  (xy+x^{u}) \cdot \adj(I) + (xy+ y^{u} ) \cdot  \adj(K) $. Hence, 
\begin{align}
\label{jointredandcore} \nonumber
     (xy+x^{u}) \cdot \adj(I) + (xy+ y^{u} ) \cdot  \adj(K) 
=&~   (xy+x^{u}, xy + y^{u})\m  &  \mbox{[by (\ref{formula for adj general})]}\\ \nonumber
=&~ (xy+ x^{u}, xy+y^{u}) (x,y) & \\  \nonumber
=&~ (x^2y + x^{u+1}, x^2y + xy^{u}, xy^2 + x^{u} y, xy^2 + y^{u+1}) & \\
=&~ ( x^{u+1},    x^2y, xy^2 ,    y^{u+1}). & 
\end{align}
We now  compute $\adj(IK)$. 

 By (\ref{eqn for IK general}),  $IK = (x^{2u+1}, x^{u+1}y, x^2y^2, xy^{u+1}, y^{2u+1})$  and hence  
\begin{align}
\label{transform-product} \nonumber
      &~ \adj (IK) R \left[  \frac{\m}{x}\right]  \\ \nonumber
=&~  \adj \left(x^{2u+1}, x^{u+2} \left( \frac{y}{x} \right), x^4 \left( \frac{y}{x}\right)^2, x^{u+2} \left( \frac{y}{x}\right)^{u+1}, x^{2u+1} \left( \frac{y}{x}\right)^{2u+1} \right) R \left[  \frac{\m}{x}\right] \\ \nonumber
=&~ x^4 \adj (K_{u}(x, y/x))R \left[  \frac{\m}{x}\right]  & \mbox{[Remark~\ref{lipman 1.2}, (\ref{defn o f Kn})]}\\
=&~ x^4 ( x^{u-2}, y/x) R \left[  \frac{\m}{x}\right].  & \mbox{[by (\ref{adj pf Kn})]}
\end{align}
By symmetry, 
\begin{align}
\label{transform-product-y}
   \adj(IK) R \left[  \frac{\m}{y}\right] 
  = y^4 ( y^{u-2}, x/y) R \left[  \frac{\m}{y}\right].
\end{align}     
Hence 
\begin{align}
\label{adjoint-product} \nonumber
   &~      \adj(IK)   & \\ \nonumber
=&~ \left ( \frac{1}{x} \adj (I K)  R \left[  \frac{\m}{x}\right]  \cap R\right) 
        \bigcap \left ( \frac{1}{y} \adj (I K)  R \left[  \frac{\m}{y}\right]  \cap R\right) &  \mbox {[by (\ref{rhodes formula})]}\\  \nonumber
=&~ \left ( \frac{x^4}{x}  \left(x^{u-2}, \frac{y}{x} \right)  R \left[  \frac{\m}{x}\right]\cap R \right)
       \bigcap \left ( \frac{y^4}{y}   \left(  \frac{x}{y}, y^{u-2} \right)   R \left[  \frac{\m}{y}\right] \cap R \right) 
     & \mbox{[by (\ref{transform-product}), (\ref{transform-product-y})]} \\ \nonumber
=&~   \left(  (x^{u+1}, x^2 y)  R \left[  \frac{\m}{x}\right] \cap R \right)
\bigcap  \left( (xy^2 , y^{u+1})      R \left[  \frac{\m}{y}\right] \cap R \right)  & \\  \nonumber
=&~ (\m^{u+1} +  y\m^2) \cap (x\m^2 +  \m^{u+1})\\ \nonumber
=&~ (x^{u+1}, x^2y, xy^2, y^3) \cap (x^3, x^2y, xy^2, y^{u+1}) & \\  
=&~ (x^{u+1}, x^2y, xy^2, y^{u+1}).
\end{align}
By (\ref{jointredandcore}) and  (\ref{adjoint-product})
\beqn
   \adj(IK) 
= (xy+x^{u}) \cdot \adj (I) + (xy+ y^{u} )\cdot \adj(K).
       \eeqn  
       We now compute $\lambda(R/  \adj(I^rK^s) )$. 
From (\ref{formula for adj general}) and Proposition~\ref{joint-reduction-zero},  for all $r,s \ge 1,$ 
\begin{align}
\label{formula for adj genreal} \nonumber
    \adj(I^rK^s) 
 =&~ I^r \cdot \adj(K^s) +\adj(I^r) \cdot K^s \\ \nonumber
= &~ I^r \m K^{s-1} + \m I^{r-1}K^s\\  \nonumber
= &~ \m I^{r-1}K^{s-1} (I+K) \\
=&~ \m I^{r-1}K^{s-1} (x^u, xy, y^u).
\end{align}
Since the ideals concerned are all complete ideals, applying Theorem~\ref{hoskin-deligne}(\ref{hoskin-deligne-b}), (\ref{multiplicity  formula general}),   (\ref{length for r general}) and from Theorem~\ref{hoskin-deligne}(\ref{hoskin-deligne-c}) we get 
\begin{align}
\label{first coefficient r I}
        e_1(I) 
= &~  e(I) - \lambda\left( \frac{R}{I}\right)  = (2u+1) - 2u = 1 \\ \label{first coefficient r I-1}
       e_1(K) 
= &~ e(K) -  \lambda\left( \frac{R}{K}\right)  = (2u+1) - 2u = 1 \\\label{mix mult of m and I general}
         e_1(\m |I) 
=&~ e_1(\m |K)=2. \hspace{.2in}
\end{align}
Put $J =(x^u, xy, y^u)$.   Then
\begin{align}
\label{length of J}
\lambda \left( \frac{R}{J} \right) = 1 + 2(u-1) = 2u-1.
\end{align} 
Similar to the computations earlier, one can verify that 
\begin{align}
\label{j red IJ}
        IJ  
= &~ (xy + y^u) I + (xy + x^u) J \\ \label{j red JK}
     JK 
= &~ (xy + x^u) K + (xy + y^u) J \\ \label{mix m IJ}
      e_1(I|J)  
=&~ e(xy + x^u, xy + y^u) = 2u\\ \label{mix mult KJ}
       e_1(K|J) 
=&~ 2u.
\end{align}
Repeatedly applying  Theorem~\ref{lipman-length-twoideals}(\ref{lipman-length-twoideals-1}) and  from (\ref{first coefficient r I}), $\cdots $, 
(\ref{mix mult KJ}) we get 
\begin{align}
\label{length formula for IrKs arbitrary}  \nonumber
&~    \lambda \left( \frac{R}{ \m I^{r-1} K^{s-1}J} \right)  & \\ \nonumber
=&~  \lambda \left( \frac{R}{\m} \right)  
+      \lambda \left( \frac{R}{I^{r-1}} \right) 
+     \lambda \left( \frac{R}{K^{s-1}} \right) 
+     \lambda \left( \frac{R}{J} \right) &\\ \nonumber
&+   e_1 (\m| I^{r-1}) + e_1 (\m | K^{s-1}) + e_1(\m|J) + e_1 (I^{r-1}| K^{s-1}) + e_1 (I^{r-1}|J) + e_1(K^{s-1}|J)\\ \nonumber
=&~ 1 + \left[  (2u+1) \binom{r+1}{2} - (2u+1) r - (r-1)\right]
 +  \left[  (2u+1) \binom{s+1}{2} - (2u+1) s - (s-1)\right]  + (2u-1)\\ \nonumber
 &+ 2(r-1) + 2(s-1) +2+  (r-1) (s-1) 2u+ 2u(r-1) + 2u(s-1)\\
 =&~ (2u+1) \binom{r+1}{2}  + 2urs +   (2u+1) \binom{s+1}{2}  -2ur -2us.
\end{align}

\noindent (d) We compute $\core (I^rK^s)$ and
$\lambda (R/\core(I^rK^s))$. 
Applying (\ref{formula for adj general}) we get
\begin{align*}
    \core (I^rK^s)  = I^rK^s \cdot \adj(I^rK^s )  = I^rK^s \m I^{r-1}K^{s-1}J  = \m I^{2r-1}K^{2s-1}J.
\end{align*}
Replacing $r-1$ by $2r-1$ and $s-1$ by $2s-1$ in  (\ref{length formula for IrKs arbitrary}) we get
\begin{align*}
\label{length for core  of  IrKs} \nonumber
   &~ \lambda \left( \frac{R}{ \m I^{2r-1} K^{2s-1} J} \right) \\ \nonumber
  =&~ (2u+1) \binom{2r+1}{2}  + 8urs +   (2u+1) \binom{2s+1}{2}  -4ur -4us\\ \nonumber
  =&~ (2u+1)  \left(4 \binom{r+1}{2}  -r \right) + 8urs + (2u+1)  \left(4 \binom{s+1}{2}  -s \right)  -4ur -4us\\
  = &~4 (2u+1)  \binom{r+1}{2} +  8urs +4  (2u+1)   \binom{s+1}{2} 
   -(6u+1) r - (6u+1)s.
        \end{align*}
}
\end{example}

  \end{document}